\def\d{\delta}

\def\hj{\hat \jmath}
\def\s3{\sqrt{3}}
\def\tr{{\mathrm tr\ }}

\def\c0{{\underline 0}}

\def\a{\alpha}

\def\Q3{{\bf Q}3} 
\def\F{{\bf F}}
\def\Q{{\bf Q}}

\def\W{{\mathcal W}}
\def\M{{\mathcal M}}
\def\O{{\mathcal O}}
\def\CS{{\mathcal S}}

\def\a{\alpha}

\def\Z{{\bf Z}} 
\def\C{{\bf C}}
\def\N{{\bf N}}

\def\lc{\lceil}
\def\rc{\rceil}
\def\lf{\lfloor}
\def\rf{\rfloor}
\def\pf{\noindent{\sc Proof. }}
\def\qed{\hfill $\Box$ \\}
\def\e{\epsilon}
\def\mod{{\ \mathrm{mod} \ }}

\documentclass[12pt]{article}
\usepackage{latexsym}
\newtheorem{dfn}{Definition}[section]
\newtheorem{thm}{Theorem}[section]
\newtheorem{lem}[thm]{Lemma}
\newtheorem{prop}[thm]{Proposition}
\newtheorem{cor}{Corollary}[thm]
\newtheorem{conj}[thm]{Conjecture}
\newtheorem{rem}{Remark}[thm]
\begin{document}
\title{Bounding slopes of $p$-adic modular forms}
\author{Lawren Smithline}
\date{June 21, 2001}
\maketitle
\abstract{
Let $p$ be prime, $N$ be a positive
integer prime to $p$, and $k$ be an integer.
Let $P_k(t)$ be the characteristic series for Atkin's $U$ operator as an
endomorphism of $p$-adic overconvergent modular forms of tame level $N$
and weight $k$.  Motivated by conjectures of Gouv\^ea and Mazur, we
strengthen a congruence in \cite{w} between coefficients of $P_k$ and
$P_{k'}$ for $k'$ $p$-adically close to $k$. For $p-1 \mid 12$, $N=1$,
$k=0$, we compute a matrix for $U$ whose entries are coefficients in the
power series of a rational function of two variables.  We apply this
computation to show for $p=3$ a parabola below the Newton polygon $\N_0$
of $P_0$, which coincides with $\N_0$ infinitely often.  As a
consequence, we find a polygonal curve {\em above} $\N_0$.
This tightest bound on $\N_0$ yields the strongest congruences
between coefficients of $P_0$ and $P_k$ for $k$ of
large $3$-adic valuation.
}

\section{Overview and background}

Let $p$ be a prime number, $N$ be a positive integer relatively prime to
$p$, and $k$ be an integer.  Let $B$ be a $p$-adic ring between $\Z_p$ and
$\O_p$, the ring of integers in $\C_p$.  Denote by $\M_k(N,B)$ the $p$-adic
overconvergent modular forms of
tame level $N$ and weight $k$ and by $\CS_k(N,B)$ the subspace of
overconvergent cusp forms.

For every weight $k$, Atkin's $U$ operator is an endomorphism of
$\M_k(N,B)$ 
stabilizing $\CS_k(N,B)$.  Denote by $U^{(k)}$ the restriction of $U$
to $\M_k(N,B)$ and by $U_{(k)}$ the restriction to $\CS_k(N,B)$.
These are compact operators, so the characteristic series
$$P_k(t) = \det(1 - tU^{(k)}), \ \ Q_k(t) = \det(1 - tU_{(k)})$$
exist.

Let $a_m(P_k)$ be the coefficient of $t^m$ in $P_k(t)$.  As a function on a
suitably defined space of weights $k$,
$a_m(P_k)$ is a rigid analytic function of $k$.

Wan\cite{w}, and Buzzard\cite{b} construct $\hat \N(m)$, which grows as
$O(m^2)$ and depends on $p$ and $N$
and not on $k$ such that $v_p(a_m(P_k)) > \hat \N(m)$.

Gouv\^ea and Mazur\cite{gm2} note, in an earlier work, the existence of
$\hat \N(m)$ and show, for prime $p \geq 5$, integer $l$ and
positive integer $n$, 
\begin{equation}\label{gmresult}
v_p(a_m(P_k) - a_m(P_{k+lp^n(p-1)})) \geq n+1.
\end{equation}
Following a remark in \cite{k}, the result in 
Equation (\ref{gmresult})
extends to $p=2, 3$.

In section 2, we show
\begin{equation}\label{stco}
v_p(a_m(P_k) - a_m(P_{k+lp^n(p-1)})) \geq \hat \N(m-2)+n+1.
\end{equation}

In section 3,
for each $p=2, 3, 5, 7, 13$, $N=1$, we construct a matrix $M$ for $U_{(0)}$
with respect to an explicit basis.  We show, for $M_{ij}$ the
entries of $M$, $$\sum_{i=1}^{\infty}\sum_{j=1}^{\infty} M_{ij}x^iy^j$$
is the power series expansion of a rational function of two variables.

In section 4, we show for $p=3$, 
$$v_p(a_m(Q_0))\geq 3{m \choose 2}+2m,$$ with equality if and only if there
is positive integer $j$ such that $m = (3^j - 1)/2$.  The secant segments
joining these vertices of the Newton polygon $\N'_0$ of $Q_0$ form a 
polygonal curve {\em above} $\N'_0$.  We find evidence in support of a
conjecture in \cite{g} on the distribution of slopes of classical modular
forms.
\\

\subsection{Motivating conjectures}

The zeros of $P_k(t)$ are reciprocals of $U^{(k)}$ eigenvalues.
For rational number $\a$,
let $d(k, \a)$ denote the number of $U^{(k)}$ eigenvalues with $p$-adic
valuation $\a$.

\begin{conj}[Gouv\^ea-Mazur]\label{gmc}
Let $k$, $l$ be integers, $n$ be a positive integer, and $\a < n$.
Then $d(k,\a) = d(k+lp^n(p-1), \a)$.
\end{conj}

Wan $\cite{w}$ uses Equation (\ref{gmresult}) and the construction of
$\hat \N(m)$ to compute a quadratic concave up function $f_{Wan}(n)$ such
that the conclusion of Conjecture \ref{gmc} holds for $\a < f_{Wan}(n)$.  

The stronger congruence in Equation (\ref{stco})
together with the method of \cite{w} shows there is quadratic $f(n)$
with quadratic term smaller than that of $f_{Wan}(n)$
such that the conclusion of Conjecture \ref{gmc} holds for $\a < f(n)$

\begin{conj}[Gouv\^ea] Let $R_k$ be the multiset of slopes with
multiplicity of classical $p$-oldforms in $\M_k(N,\Z_p)$.
The probability that an element of $R_k$ chosen with uniform distribution
is in the interval $(\frac{k-1}{p+1},\frac{p(k-1)}{p+1})$ 
diminishes to zero as $k$ increases without bound.
\end{conj}


\subsection{Spaces of overconvergent modular forms}\label{allpr}

For $p \geq 5$, let $E_{p-1}$ be the level one Eisenstein series.
Let $M_k(N,B)$ be the classical weight $k$ level $N$ modular forms with
coefficients in $B$.

\begin{prop}[Katz]
For $p \geq 5$ and any $f \in \M_k(N,B)$,
there are $b_j \in M_{k+j(p-1)}(N,B)$ for $j\geq 0$ and $r \in \O_p$
of positive valuation such that 
\begin{equation}\label{katzexp}
f = b_0 +\sum_{j=1}^\infty r^j b_j / E^j_{p-1},
\end{equation}

There is a distinguished choice of $b_j$
after choosing $r$ and direct sum decompositions
$$M_{k+j(p-1)} = E_{p-1}\cdot M_{k +(j-1)(p-1)} \oplus W_{k+j(p-1)},$$
such that $b_j \in W_{k+j(p-1)}(N,B)$ for $j >0$.
\end{prop}

See \cite{k}, Propositions 2.6.1 and 2.8.1.
The parameter $r$ is 
the {\em growth condition} and $v_p(r)$ is bounded above by the given $f$.

Let $\M_k(N,B,r)$ be the space of modular forms with growth condition $r$.
The space $\M_k(N,B)$ is $\bigcup_{v_p(r) > 0} \M_k(N,B,r)$.


\begin{rem} For $p=3$, $N > 2$ and prime to 3,
Theorem 1.7.1 of \cite{k} shows there is a
level $N$ lift of the characteristic 3 Hasse invariant, so an analogous
expansion result holds.
Proposition 2.8.2 of {\em loc. cit.} shows the expansion result for $N=2$.
\end{rem}

\begin{prop} Suppose $p = 2$ or $3$ and $N$ relatively prime
to $p$. For any  $f \in \M_k(N,B)$ there are $b_j \in M_{k+4j}(N,B)$ and $r
\in \O_p$ of positive valuation such that
$$f = b_0 + \sum_{j=1}^\infty r^{4j/(p-1)} b_j/ E^j_4,$$
There is a distinguished choice of $b_j$ after choosing $r$ and direct sum
decompositions $$M_{k+4j} = E_4\cdot M_{k +4(j-1)} \oplus W_{k+4j}(N,B),$$
such that $b_j \in W_{k+4j}(N,B)$ for $j>0$.
\end{prop}

\pf
We follow the remark at the end of Subsection 2.1 of {\em loc. cit.}.
Let $B$ be the fourth power of the Hasse invariant $A$ for $p=2$ and the square
of $A$ for $p=3$.  In either case, $B$ is a weight 4 level
1 modular form defined over $\F_p$.  A version of Deligne's congruence
holds: $B \equiv E_4 \mod 2^4$ and $B \equiv E_4 \mod 3$.

For $N > 2$, and relatively prime to $p$, the functor ``isomorphism classes
of elliptic curves with level $N$ structure'' is representable by a scheme
which is smooth over $\Z[\frac1N]$ and the formation of modular forms
commutes with base change to a ring in which $p$ is topologically nilpotent.
So we repeat the construction of $p$-adic modular forms for $p=2, 3$
and Katz expansions with powers of $r^{4/(p-1)}E_4^{-1}$.

For $p=2, 3$ (and 5), and $N = 1$,
Section 1.4 of \cite{s2} states weight zero forms
have expansions in powers of $\Delta E_4^{-3}$ where $\Delta$ is the weight 12
level 1 cusp form. Coleman\cite{c} shows $$E_k \cdot\M_0(N,B) = \M_k(N,B).$$

$M_k(N,B)$ is a free $B$ module, so $M_k(N,B) = E_4 \cdot M_{k-4}(N,B)
\oplus W_k(N,B)$ for some $W_k(N,B) \subset M_k(N,B)$.
\qed

\begin{thm}[Coleman]\label{cow} Let $k_1, k_2$ be weights.
Let $G(q) \in M_{k_1 - k_2}(N,B)$.  Let $\Xi$ be the operator
multiplication by $G(q) / G(q^p)$.  
If
$1/G \in \M_{k_2-k_1}(N,B)$
then $U^{(k_1)}$ is similar to $U^{(k_2)} \Xi$.
\end{thm}

\begin{rem} 
The Eisenstein series satisfy the hypothesis of Theorem \ref{cow}.
\end{rem}

\subsection{Notations for matrices and Newton Polygons}

Let $M$ be a matrix over a ring, possibly of infinite rank.
Let $n$ be a nonnegative integer.
Let $s = (s_1, s_2, s_3, ... s_n)$ be a sequence of $n$ distinct
natural numbers.

The {\em $n \times n$ diagonal major} of $M$ associated to
$s$ is the $n \times n$ matrix $A$ whose entry $A_{ij}$ is $M_{s_i, s_j}$.

A {\em selection} of a $M$ associated to $s$ and degree $n$
permutation $\pi$ is a sequence of $n$ elements, $( M_{s_1, s_{\pi(1)} },
M_{s_2, s_{\pi(2)}}, \ldots, M_{s_n, s_{\pi(n)}}).$

The {\em $n \times n$ diagonal minor} of $M$ associated to $s$
is the determinant of the $n \times n$ diagonal major of $M$ associated to
$s$.

The {\em upper $n \times n$ diagonal major} of $M$ is the
diagonal major associated to the sequence $(1, 2, 3, \ldots n)$.

The diagonal matrix $D= diag(d_i: i\geq 1)$ is the matrix with entries
 $D_{ii} = d_i$ and zero elsewhere.

The Newton polygon of power series $P(t)$ is the function 
$\N(m)$ which is the lower convex hull of the set $(m, v_p(a_m(P)))$, defined
for real $m \geq 0$.

A vertex of the Newton polygon $\N(m)$ is a point $(m, \N(m))$
such that $\N(m) = v_p(a_m(P))$.

A side of a Newton polygon $\N(m)$ is a line segment whose
endpoints are vertices.

The slopes of a Newton polygon are the slopes of its sides.

The multiplicity of a slope is the difference of the first
coordinates of its endpoints.

We denote by $\N_k(m)$ the Newton polygon of $P_k$, and by $\hat \N_k(m)$ 
a function such that $\N_k(m) \geq \hat \N_k(m)$.  We indicate by $\hat
\N(m)$ a function such that for all weights $k$, $\N_k(m) \geq \hat\N(m)$.

We denote by $\N'_k(m)$ the Newton polygon of $Q_k$,
and by $\hat \N'_k(m)$ a function such that $\N'_k(m) \geq \hat \N'_k(m)$.

We state as Lemma \ref{genii} that if $p-1 \mid 12$ and $N=1$, then
$P_k(t) = (1-t)Q_k(t)$.  For these cases, $\N'_k(m) = \N_k(m+1)$.

\section{Comparing Newton polygons for $U$ in different weights}
Retain $p, N, k$ as before, and let $l$ be an integer and $n$ be a positive
integer.
For $p=2$, we require $n \geq 2$.
Let $k' = k + l(p-1)p^n$. At the end of the section, we show
there is a quadratic $\hat \N(m)$ such that
$$
v_p(a_m(P_k) - a_m(P_{k'})) \geq \hat \N(m-2)+n+1.
$$

We now describe only the case $p> 3$ for clarity.
Section \ref{allpr} reviews the differences for $p=2, 3$
from the case $p > 3$.

Let $r = p^{1/(p+1)}$.
Choose a basis $\{ b_{0,s} \}$ for the module $M_k(N,B)$.
For $i>0,$ choose a basis $b_{i,s}$ for the module $W_{k+i(p-1)}(N,B)$.

Let $e_{i,s} = r^i E_{p-1}^{-i} b_{i,s}$.  
Let $M$ be the matrix for
$U^{(k)}$ with respect to the basis $\{ e_{i, s} \}$.

Let $\N_k(m)$ be the Newton polygon of $P_k(t)$.

Lemma 3.1 of \cite{w} includes
\begin{lem}
Let $M_{i,s}^{u,v}$ be the coefficient of $e_{u,v}$ in $U^{(k)} (e_{i,s})$.

Then $v_p(M_{i,s}^{u,v}) \geq u (p-1)/(p+1)$.
\end{lem}

Let $d_u = \dim M_{k+u(p-1)}(N,B) \otimes \C_p$.
For $u > 0$, let $m_u = d_u - d_{u-1}$.

\begin{lem}[Wan] 
Let $k$ be a weight.  If $d_v \leq m < d_{v+1}$
for some $v\geq 0$, then
\begin{equation}\label{wanest}
\N_k(m) \geq \frac{p-1}{p+1}\left(\sum_{u=0}^v um_u+ (v+1)(m-d_v)\right) -m.
\end{equation}
\end{lem}

\begin{dfn} Let $\hat \N_k(m)$ be the right side of Equation (\ref{wanest}).
\end{dfn}

The $m_u$ have an upper bound, depending on $p$ and $N$, so
$\hat \N_k(m)$ grows quadratically.

Wan shows   $\N_k(m) = \N_{k'}(m)$ when both are less than $n+1$.








\begin{lem}\label{sizeps}
Let $A$ be the matrix for $U^{(k)}$ with respect to basis $\e_{i,s} =
r^{-i}e_{i,s}$.
Then
$$v_p(A_{i,s}^{u,v}) \geq (up -i )/(p+1)$$ and also at least zero.
\end{lem}
\pf $U^{(k)}$ stabilizes $\M_k(N,B,1)$, as shown in \cite{gm2}.
\qed

\begin{prop}\label{goodint}
$E_{p-1}^{p^n}(q) / E_{p-1}^{p^n}(q^p) \in 1+ p^{n+1}\M_0(1,\Z_p,1).$
\end{prop}
\pf In weight zero, the only $\e_{i,s}$ not 0 at the cusp
$\underline\infty$ is the constant function 1.  The $q$-expansion of
$(E_{p-1}-1)/p $ is in $q\Z[[q]]$.
\qed

\begin{thm}\label{highn}
For $k$, $k'$ as above, \ 
$v_p(a_m(P_k) - a_m(P_{k'})) \geq \hat \N_k(m-2) + n + 1.$
\end{thm}

\pf
Let $C$ be the matrix with respect to the basis $\e_{i,s}$
for multiplication by $E_{p-1}^{p^n}(q) /
E_{p-1}^{p^n}(q^p)$ considered as an operator on $\M_k(N,B,r)$.

Let $M^{(k')} = MC$.  By Theorem \ref{cow}, $M^{(k')}$ is a matrix for an operator similar to
$U^{(k')}$ on $\M_{k'}(N,B,r)$ and $M^{(k')}$  acts on $\M_k(N,B,r)$.

By Proposition \ref{goodint}, the matrix $C-1$ is a matrix with
entries in $p^{n+1}B$, so $M - M^{(k')}$ has entries in $p^{n+1}B$.

The difference $a_m(P_k) - a_m(P_{k'})$ is equal to
$$\tr \bigwedge^m M - \tr \bigwedge^m M^{(k')}.$$

These traces are the sums of all the different $m \times m$ diagonal minors
of $M$ and $M^{(k')}$, so the difference contains terms (up to sign)
\begin{equation}\label{selprod}
\prod_{i=1}^m M^{(k')}_{s_i,s_{\pi(i)}} - \prod_{i=1}^m M_{s_i,s_{\pi(i)}},
\end{equation}
where $s$ is a sequence of $m$ integers, $\pi$ is a permutation of degree
$m$.

Let
\begin{equation}\label{insprod}
Z = \prod_{i=1}^m (z_i+w_i) - \prod_{i=1}^m (z_i),
\end{equation}
where $z_i \in B$ and $w_i \in p^{n+1}B$,
be instance of equation (\ref{selprod}).

The sequence $(z_1, z_2, ... z_m)$ is a selection of $M$. By Lemma
\ref{sizeps}, the product of any $m-j$ of them has valuation at least
$\hat \N_k(m-2j)$.  The product of any $j$ of the $w_i$ has valuation at least 
$j(n+1)$.

Rewrite (\ref{insprod}) as
\begin{equation}\label{insdet}
Z = \sum_{\emptyset \neq s \subset \{1, 2, ... m\}  } \prod_{i\in s}w_i
\prod_{i \not\in s} z_i.
\end{equation}

For any subset $s$ of size $j$,
$$v_p(\prod_{i\in s} w_i \prod_{i\notin s} z_i) \geq \hat \N_k(m-2j) + j(n+1).$$

The set $s$ is nonempty, so,
$$v_p(Z) \geq \hat \N_k(m-2) + n+1,$$
for every instance of Equation (\ref{insprod}).


\qed

\begin{cor} There is a quadratic $\hat \N(m)$ independent of $k$
such that the conclusion of Theorem \ref{highn} holds.
\end{cor}

\pf
Given $p, N$, Wan\cite{w} shows
there are finitely many different $\hat \N_k(m)$.
Let $\hat \N(m)$ be the infimum of them.
\qed

\section{Computing tame level 1 $U$ for $p\in\{2,3,5,7,13\}$}
\label{comu}

Let $p$ be a prime such that $X_0(p)$ has genus 0, that is, $p \in \{2, 3,
5, 7, 13 \}$ and $N=1$.  We show how to compute $U_{(0)}$ with respect to
an explicit basis.

The curve
$X_0(p)$ has a uniformizer
$$d_p = \sqrt[p-1]{\Delta(q^p)/\Delta(q)}$$
with simple zero at the cusp $\underline \infty$, pole at the cusp
$\underline 0$, and leading $q$ expansion coefficient $1$.  

Let $\pi:X_0(p) \rightarrow X_0(1)$ be the map which ignores level $p$
structure.  Let $\hj = \pi^*(j)$.
The map $\pi$ is ramified above $j =
0, 1728, \infty$ only.  

\begin{prop}\label{polex}
There is a degree $p+1$ polynomial $H_p$ over $\Z$
with constant term 1 such that
$$d_p \hj = H_p(d_p).$$
\end{prop}
\pf
The map $\pi$ has degree $p+1$. The product $d_p \hj$ has a pole
only at the cusp $\underline 0$. Hence, there is a polynomial $H_p$
satisfying the proposition.

$H_p$ has integer coefficients, because
the $q$-expansion of $d_p\hj$ at $\underline\infty$
is in $1+q\Z[[q]]$
\qed

\begin{rem}
The ramification degrees of $\pi$ over $j=0$ are $1$ and $3$,
yielding roots of multiplicity 1 or 3 of $H_p(d_p)$.  Points over $j=1728$ are
roots of multiplicity 1 or 2 of $H_p(d_p) - 1728d_p$.  We calculate
$H_p$ by equating $q$-expansions.
\end{rem}


\begin{lem}\label{genii} $P_k(t) = (1-t)Q_k(t)$
\end{lem}
\pf
$X_0(p)$ has genus 0, so the only weight zero noncuspidal eigenforms
are constants and the eigenvalue is 1.  By a theorem of \cite{h},
or as a consequence of Theorem \ref{cow}, in every weight $k$, $d(k,0)=1$
and a slope zero eigenform is noncuspidal.
\qed

Let $t_2=4$, $t_3=3$.  For $p \geq 5$, let $t_p = 1$.

Let $c_2 =0$, $c_3=1728$, $ c_5 = 0$, $c_7=1728$, and $c_{13} = 432000/691$.

Let $e=12/(p^2-1)$.  

\begin{lem}\label{newhook}
The Newton polygon of $H_p(d_p) - c_pd_p$, as a polynomial in $d_p$,
 has a single side of slope
$ep$.
\end{lem}

\begin{lem}\label{poweis}
The weight 12 power of $E_{t_p(p-1)}$ is $(j - c_p)\Delta.$
\end{lem}

The lemmas are direct computations.

\begin{prop}
Let $r < p/(p+1)$.  The disc $D = \{z: z \in X_0(1), v_p(E_{t_p(p-1)}(z)) <
t_p r \}$
is isomorphic to $\{z: z\in X_0(p), v_p(d_p(z)) > -er(p+1)\}$.
\end{prop}
\pf
When $z \in X_0(1)$ is a point of supersingular reduction, $\Delta(z)$ is a
unit.  At a point of ordinary reduction, $E_{t_p(p-1)}(z)$ is a unit and
$v_p(\Delta(z) ) \geq 0$.  By Lemma \ref{poweis},
$D = \{z: v_p(j(z) - c_p) < er(p+1)\}$.

Lemma \ref{newhook} shows the relation
$(\hj-c_p)d_p = H_p - c_pd_p$
is uniquely invertible for $d_p$
such that $v_p(d_p(z)) > -er(p+1)$, establishing the isomorphism.
\qed

\begin{cor}  $\CS_0(1,\Z_p) \subset d_p\Z_p[[d_p]].$
$U_{(0)}$ acts as a matrix $M$ on a basis of powers of $d_p$.
\end{cor}

Let $\W$ be the rigid subspace of $X_0(p)$ where $v_p(\pi^*(E_{t_p(p-1)})) <
t_p/(p+1)$. 
The section $s$ of $\pi$ over $\pi(\W)$ such that for elliptic curve
$E$, $s(E)$ is the pair $(E, C)$ for $C$ the canonical order $p$ subgroup
of $E$ is an isomorphism.  

Let $V$ be the pullback of $\phi$, the Deligne-Tate lift of
Frobenius on $X_0(1) / \F_p$.
Let $w_p$ be the Atkin-Lehner involution on $X_0(p)$.

\begin{lem}
For points of $\W$, 
$$ V( j) \circ \pi = \hj \circ w_p.$$
\end{lem}
\pf
The Atkin-Lehner involution acts as
$$w_p :(E,C) \rightarrow (E/C, E[p] / C).$$

$E$ has a canonical subgroup of order $p$, and
$$V: E \rightarrow E/ \ker \phi^*,$$
coincides with $s^*\circ w_p^* \circ \pi^* $.
\qed
We identify $\W$ with $\pi(W)$ via section $s$.

\begin{prop}
For points of $\W$,
\begin{equation}\label{modeq}
H_p(p^{12/(1-p)}/ d_p) V(d_p) - p^{12/(1-p)}H_p(V(d_p)/d_p = 0.
\end{equation}
\end{prop}

\pf The modular equation
$$H_p(w_p^*(d_p)) V(d_p) = H_p(V(d_p))w_p^*(d_p)$$
holds on $\W$ and $w_p(d_p) = (p^{12/(1-p)} / d_p)$.
\qed

\begin{thm}\label{algfn}
There is an algebraic function $I_p(y,x)$ and a matrix $M$ for $U_{(0)}$
with respect to the basis $d_p^n$ such that entries $M_{ij}$
 satisfy a generating function equation
\begin{equation}\label{mgeneq}
\sum_{i = 1}^{\infty}\sum_{j=1}^{\infty} M_{ij}x^iy^j =
\frac{y}{p} \frac{d}{dy} \log I_p(x,y).
\end{equation}
\end{thm}

\pf
Clear denominators and factor $V(d_p) - w_p^*(d_p)$ from Equation
(\ref{modeq})
to determine an algebraic relation  $$d_p^p I_p(V(d_p), 1/d_p)$$
between $d_p$ and $V(d_p)$, 
of degree $p$ in $d_p$.  The inverse of $V$ applied to coefficient of
$d_p^{p-1}$ is $\tr V(d_p) = p U(d_p)$.

The values of $U(d_p^n)$ for $n = 0$ to $p-1$ and the coefficients of $I_p$
determine a recurrence for $U(d_p^n)$ for $n \geq p$.  
\qed

\begin{rem}\label{theip}
The $I_p(x,y)$ for $p=2,3,5,7, 13$ are \end{rem}
\begin{eqnarray*}
I_2 & = & 1 - (2^{12}x^2 + 3\cdot 2^4x)y - xy^2, \\
I_3 & = & 1 - (3^{12}x^3 + 4\cdot 3^8x^2 + 10\cdot 3^3x)y
            - (3^6x^2 + 4\cdot3^2x)y^2 - xy^3, \\
I_5 & = & 1 - (5^{12}x^5 + 6\cdot 5^{10}x^4 + 63\cdot5^7x^3 + 52\cdot5^5x^2
+ 63\cdot5^2x)y \\ 
& & - (5^9x^4 + 6\cdot5^7x^3 + 63\cdot5^4x^2 + 52\cdot5^2x)y^2 \\
& & - (5^6x^3 + 6\cdot5^4x^2 + 63\cdot5x)y^3 - (5^3x^2 + 6\cdot5x)y^4 -
xy^5, \\
I_7 & = & 1 - (7^{12}x^7 + 4\cdot7^{11}x^6 + 46\cdot 7^9x^5 + 
272\cdot7^7x^4 + \\ & &
845\cdot7^5x^3 + 176\cdot7^2x^2 + 82\cdot7x)y - ... - xy^7, \\
I_{13} & = & 1 - (13^{12}x^{13} + 2\cdot13^{12}x^{12} +
25\cdot13^{11}x^{11} + 196\cdot13^{10}x^{10} + \\
& &  1064\cdot13^9x^9 + 4180\cdot13^8x^8 +
12086\cdot13^7x^7 + \\
& & 25660\cdot13^6x^6 + 39182\cdot13^5x^5 +  
41140\cdot13^4x^4 + \\
& & 27272\cdot13^3x^3 + 9604\cdot13^2x^2 + 1165\cdot13x)y - ... - xy^{13}. 
\end{eqnarray*}



\begin{prop}\label{simpar}
The $p$-adic valuation of $M_{ij}$ is at least $e(pi - j)-1$.
There is a parabola $\hat \N(m)$
with quadratic coefficient $6/(p+1)$ such that $\N_0(m) \geq \hat \N(m)$.
\end{prop}
\pf
Let $M'_{ij} = p^{e(j-i)}M_{ij}.$  The matrix $(M'_{ij})$ is similar
to $(M_{ij})$.
Theorem \ref{algfn} shows
\begin{equation}\label{yacc}
\sum_{i=1}^\infty \sum_{j=1}^\infty
M'_{ij}x^iy^j = \frac{y}{p}\frac{d}{dy} \log I_p(p^{-e}x,p^ey).
\end{equation}
Direct calculation shows $I_p(p^{-e}x,p^ey)$, for the $I_p$
displayed in Remark \ref{theip}, is a polynomial
in $p^{e(p-1)}x$ and $y$ with integer coefficients.
Hence, $v_p(M'_{ij}) \geq i\cdot e(p-1).$
\qed

\subsection{Tame level 1 and $p=2$ or $3$}

Emerton\cite{e} calculates the lowest positive slope $2$-adic modular
forms of every weight.  Concise expressions for the $q$-expansions of a few
forms facilitate computation.

Serre\cite{s} observes that for a compact operator $M$ expressed as a
matrix on a basis of a Banach space, if $c_i$ is the infimum of the
valuations of column $i$ of $M$, then $\tr(\wedge^n M)$ has valuation
at least the sum of the $n$ smallest $c_i$.

  
\begin{prop} For $p=2$ and even weight k, 
there is an $\O_2$ basis $\{e_n\}_{n\geq 1}$ of $\CS_k(1,
\O_2)$ such that the image of
$U_{(k)}$ is a subset of $\bigoplus 8^n e_n \O_2.$
\end{prop}
\pf  This is a rewriting of Proposition 3.21 of \cite{e} in language
amenable to the noted observation of Serre.  The basis element $e_n$ is
$F_k d_2^n$ for a certain weight $k$ form $F$.
\qed

Recall $\N'_k(m)$ is the Newton polygon of $Q_k(t)$.

\begin{cor} $\N'_k(m) \geq 3{m+1 \choose 2}$.
\end{cor}

\begin{lem}\label{whs}
Suppose $p=3$.
Let $S = \sqrt[8]{\Delta^3 / V(\Delta)}$.
$S^2$ is in $\M_6(1, \Z_p)$ and does not vanish at the cusp
$\underline\infty$.
The quotient $S / V(S)$ is in $M_0(1, \O_3, 3/2)$ and as a power series in
$Z[[d_3]]$, $S/V(S) -1$ is in the ideal $(9d_3, 27d_3^2)$.
\end{lem}
\pf
Direct calculation and comparison of $q$ expansions shows $S$ is the
Eisenstein series for level 3, weight $3$ and character $\tau$, the $3$-adic
Teichmuller character.  $S^2$ is a level 3 weight 6 classical modular
form and thus a tame level 1 weight 6 overconvergent modular form.

The curve $X_0(9)$ has genus zero and uniformizer 
$$d_9 = \sqrt[8]{V(V(\Delta))/\Delta}.$$
The ramification of the forgetful map to $X_0(3)$ shows
$$d_3 = d_9 + 9d_9^2 + 27d_9^3.$$

Reversal of this relation between $d_3$ and $d_9$ and the observation
$$S/V(S) = d_9 / d_3$$  shows $S/V(S)$ is in $M_0(1,\O_3, 3/2)$,
has constant term 1, and $S/V(S) -1 \in (9d_3, 27d_3) \subset Z[[d_3]]$.
\qed

\begin{prop}  For $p=3$ and even weight $k$ divisible by 3,
$\N'_k(m) \geq 3{m \choose 2}$.
\end{prop}
\pf
Let $R$ be the multiplication by $(S/V(S))^{k/3}$ operator.
Theorem \ref{cow} shows the composition $U_{(0)}R$
is similar to $U_{(k)}$.
Lemma \ref{whs} shows the conclusions of Proposition \ref{simpar}
hold for $U_{(0)}R$.
\qed

\subsection{Further example for $p=3$, $N=1, k=0$.}

Let $p=3$, $N=1$ and
$$\hat \N'_0(m) = \frac32 m (m-1) + 2m.$$
We work an example of Proposition \ref{simpar}.

\begin{lem}\label{furex}
$\N'_0(m) \geq \hat \N'_0(m)$.
\end{lem}

\pf
Recall $e = 3/2.$  Equation (\ref{yacc}) shows
\begin{equation}\label{ngeneq}
3\sum_{i,j} M'_{ij}x^iy^j =
\frac{9(10xy + 8\sqrt{3}xy^2 + 3xy^3) +
3^5(4\sqrt{3}x^2y + 2x^2y^2) + 3^8x^3y}
{1- 3^3(10xy + 4\sqrt{3}xy^2 + xy^3)
- 3^6(4\sqrt{3}x^2y + x^2y^2) -3^9x^3y}.
\end{equation}
Following the last step of Proposition \ref{simpar},
substitute $\d = 3^3x$ into the right side of
Equation (\ref{ngeneq}) to get
\begin{equation}\label{gdef}
G(\d,y) =  \frac{10\d y + 8\sqrt{3}\d y^2 + 3\d y^3 +
4 \sqrt{3}\d^2y + 2\d^2y^2 + \d^3y}
{1 - 10\d y - 4\sqrt{3}(\d y^2 + \d^2y) - (\d y^3+\d^2y^2+ \d^3y)}.
\end{equation}

The valuation of $M'	_{ij}$ is at least $i \cdot e(p-1)-1 = 3i-1$.
So $$\N'_0(m) \geq \sum_{i=1}^m 3i-1 =\hat \N'_0(m).$$
\qed

\section{For $p=3$, $N=1$, $\hat \N'_0$ is a sharp parabola below $\N'_0$}

Let $p =3$ and $N=1$ and
$$m_i = \sum_{j=0}^{i-1} 3^j = \frac{3^i-1}2.$$

\begin{thm} 
The set $E = \{m: m \in \Z, \N'_0(m) = \hat \N'_0(m) \}$
is the same as $\{ m_i : i \geq 0 \}.$
\end{thm}

\pf  We show for all $m \geq 0$, that
$m \in E$ if and only if $(m - 1)/3 \in E$.

The leading coefficient of $P_0$ is 1, so $0 \in E$.
\\

Let $M'$ be the matrix for $U_{(0)}$ with respect to basis $\{ 3^{3m/2}d^m
\}$.

Lemma \ref{furex} shows $M'_{ij}$ has valuation at least $3i-1$, so
there is a matrix $K$ over $Z[\sqrt{3}]$
and diagonal matrix $D = diag(3^{3i-1})$  such that $M' = DK$.

Let $\bar K = K \mod \sqrt{3}\Z[\sqrt{3}]$ and let
$c_m(\bar K)$ be its upper $m \times m$ diagonal minor.


Every $m \times m$ diagonal minor of $M'$ has valuation at least $\hat
\N'_0(m)$ and the inquality is strict except for the upper $m \times m$
diagonal minor.  So we have reduced the theorem to showing
that $m \in E$ if and only if $c_m(\bar K) \neq 0$.
\\

Call a degree $m$ permutation $\pi$ {\em excellent} if the selection of $\bar
K$ associated to $(1, 2, \ldots m)$ and $\pi$ is a sequence of nonzero
entries of $\bar K$.

\noindent{\bf Claim 1.}  If there is a degree $m$ excellent $\pi$, then
$m = m_i$ for some $i$.

We establish Claim 1 by induction.  The trivial degree 0 permutation is
excellent.

The entries of
$K$ satisfy a linear recurrence.
Equation (\ref{gdef}) with $x$ substituted for $\delta$ is
$$G(x,y)= 
\frac{10xy + 4\s3 xy(x+2y) + xy(x^2 + 2xy + 3y^2)}{1-xy(10+4\s3(x+y) +
x^2+xy+y^2)}.$$
The coefficient of $x^iy^j$ is the entry of $K$ in row $i$ and column $j$.

Let $\bar G$ be the generating function for entries of $\bar K$.
$\bar G$ is the reduction of $G$ to $\F_3[[x,y]]$.

Let  
$$R(i) = (1+(xy+x^3y+x^2y^2+xy^3)+(xy+x^3y+x^2y^2+xy^3)^2)^{3^i},$$
and
$$\bar G_0(x,y) = xy(1-xy + y^2).$$

Let $$\bar G_j = \bar G_0 \cdot \prod_{i=0}^{j-1} R(i)$$
 and $$\bar C_j =
\prod_{i=j}^\infty R(i).$$

For all nonnegative integers $j$, $\bar C_j^3 = \bar C_{j+1}$ and
$\bar G = \bar G_j \bar C_j$.

By direct computation,
\begin{equation}\label{recg}
\bar G_1 = (x^{-1}y+1 - xy^{-1} + y^{-2}) \bar G_0^3
+ xy+x^2y^4+x^6y^2,
\end{equation}
and so
\begin{equation}\label{gwrite}
\bar G = (x^{-1}y + 1 - xy^{-1} + y^{-2}) \bar G^3 + (xy + x^2y^4 +
x^6y^2)\bar C_1.
\end{equation}

Equation (\ref{gwrite}) shows the coefficient of $x^iy^{3j}$ in $\bar G$ 
is the same as the coefficient of $x^iy^{3j}$ in $\bar G^3$.  This
coefficient is zero if $i$ is not divisible by 3.

Suppose degree $m$ permutation $\pi$ is excellent.
The only unit in row 1 is in column 1, so $\pi(1)=1$.
The functions
\begin{equation}\label{splitpi}
\sigma(i) = \pi(3i)/3, \ \ \sigma'(i) = (\pi(3i-1) -1)/3, \ \ 
\sigma''(i) = (\pi(3i+1)+1)/3
\end{equation}
 are
excellent degree $\lf m/3 \rf$ permutations, and $3 \mid (m-1)$.

The inductive step is complete.

\noindent
{\bf Claim  2.} For any $m_i$, there is a unique degree $m_i$
excellent $\pi$.

We proceed by induction.  The unique degree 0 permutation is excellent.

Equation (\ref{gwrite}) shows for excellent degree $\frac{m-1}3$
permutations $\sigma, \sigma', \sigma''$, there is an excellent degree
$m$ permutation $\pi$, computed by reversing Equations (\ref{splitpi}).

If there is a unique degree $(m-1)/3$ excellent $\sigma$,
then there is a unique degree $m$ excellent $\pi$.
Claim 2 is established.

Claim 1 shows for $m$ not equal to any $m_i$, that $c_m(\bar K) = 0$.
Claim 2 shows for each $m_i$, there is a unique selection of the upper $m_i
\times m_i$ diagonal major of $\bar K$ which contributes a nonzero term to
$c_{m_i}(\bar K)$.  Hence, $c_m(\bar K) \neq 0$ if and only if there is $i$
such that $m = m_i$.
\qed

\begin{cor}\label{above}
Let $L$ be the secant line such that $L(m_i)= \hat \N'_0(m_i)$ and
$L(m_{i+1})= \hat \N'_0(m_{i+1})$.  If  $m$ is such that $m_i < m < m_{i+1}$,
then $$\hat \N'_0(m) < \N'_0(m) \leq L(m).$$
\end{cor}



\begin{prop}\label{binob}
Let $l$ be an integer, $n$ be a nonnegative integer.
Let $k= 2\cdot 3^{n+1} \cdot l$.
Let $s$ be an integer, $0 \leq s < 2\cdot3^{n-1}.$
If $\N'_0(s) = \hat \N'_0(s)$,
then $\N'_k(s) = \hat \N'_0(s).$
\end{prop}
\pf
Let $R = (S/V(S))^{k/3}$.  The binomial theorem shows
the coefficient of $d_3^m$ in $R$
has valuation at least $\lc 3m/2 \rc + n - v_3(m)$.

Let $C$ be the matrix for the
multiplication by $R$ operator on $\CS_0(1, \O_p)$ 
with respect to the basis $\{ 3^{3m/2}d^m \}$.

Let $M'$ be the matrix for $U_{(0)}$ with respect to the 
same basis.  

By Theorem \ref{cow}, $M'C$ is similar to a matrix for $U_{(k)}$.

For all $i, j$, $v_3(M'_{ij}) \geq 3i-1$.
For $i > 3j$ or $j > 3i$, $M'_{ij} = 0$.

For all $j>0$, $C_{jj} = 1$.
For $j, m>0$, $v_3(C_{j+m, j}) \geq n - v_3(m)$ and $C_{j,j+m}=0$.

For odd $m$, including $m=3^{n-1}$, $v_3(C_{j+m , j}) \geq \frac12$.

For all $i$,  $v_3(M'_{ij} -(M'C)_{ij}) \geq 3i-1.$

For $i \leq s$, $v_3(M'_{ij} - (M'C)_{ij}) > 3i-1$,
because $$(M'C)_{ij} = \sum_{k=j}^{3i} M'_{ik}C_{kj},$$
and $3i \leq 3s < 2\cdot 3^n$.

If $\N'_0(s) = \hat \N'_0(s)$ then $\N'_k(s) = \N'_0(s)$.
\qed



\begin{cor}
Let $l$ be an integer and $n$ be a nonnegative integer.  Let $k=2\cdot
3^{n+1}\cdot l$.

For integer $i$, $0 \leq i <
n-1$, there are exactly $3^i$ overconvergent 3-adic modular forms of
weight $k$ with slope in $[m_{i+1}+1, m_{i+2}-2]$, and these have
average slope $3^{i+1} - 1$.
\end{cor}

\pf By Proposition \ref{binob}, $\N'_k(m_i) = \N'_0(m_i)$ and
$\N'_k(m_{i+1}) = \N'_0(m_{i+1})$.
There are $3^i = m_{i+1} - m_i$ slopes with multiplicity accounted for by
the edges joining these vertices of the Newton polygon $\N'_k$. 
The difference $\N'_k(m_{i+1}) - \N'_k(m_i)$
is $3^i(3^{i+1} -1)$.

The average slope is $3^{i+1} -1$.  The minimum of these slopes is
at least $3m_i+2$ and the maximum at most $3m_{i+1}-1.$ \qed

\begin{cor}
Let $k$ be an even integer and $i$ be a positive integer. If 
$$v_3(k) \geq [\hat \N'_0(m_{i+1}) + \hat \N'_0(m_{i})]/2 - 
\hat \N'_0((m_{i+1}+m_{i})/2) + i+2,$$
then for $m \leq m_{i+1}$, $\N'_0(m) = \N'_k(m)$.
\end{cor}
\pf  
The Newton polygons $\N'_0$ and $\N'_k$ both have vertices $(m_i, \hat \N'_0(m_i))$
and $(m_{i+1}, \hat \N'_0(m_{i+1}))$.


By Corollary \ref{above} and Theorem \ref{highn},
$v_3(a_m(P_k)) = v_3(a_m(P_0))$ for every $m$ between $m_i$ and $m_{i+1}$.
\qed

Affirming a pattern noticed by Gouv\^ea\cite{g},
\begin{cor} Let $k = 2\cdot 3^{n+1}$.
The classical weight $k$ level 3 oldforms have slopes outside
$[k/4,3k/4]$.
\end{cor}
\pf  There are $m_n = \frac{k}{12} - \frac12$ cuspidal level 1 normalized
eigenforms.  There are $2m_i+2$ classical level 3 oldforms, and one pair of
these comes from the weight $k$ Eisenstein series.  The slopes of the forms
in this pair are $0$ and $k-1$.

By Proposition \ref{binob}, $\N'_k(m_n) = \hat \N'_0(m_n)$, because $m_n <
2\cdot 3^{n-1}.$

The slope $\N'_k(m_n) - \N'_k(m_n -1)$ is less than
$$\hat \N'_0(m_n) - \hat \N'_0(m_n-1) = 3m_n - 1 = \frac{k}4 -1.$$

The mates of these $m_i$ oldforms have slopes greater than $\frac{3k}4$.
\qed

\noindent
Department of Mathematics \\
Malott Hall \\
Cornell University \\
Ithaca, NY 14853 \\
USA \\
\\
{\tt lawren@math.cornell.edu}

\end{document}